\documentclass[12pt]{iopart}

\usepackage{iopams,color}
\usepackage{graphicx,cite}

\newtheorem{thm}{Theorem}[section]

\newtheorem{lem}[thm]{Lemma}

\newtheorem{rek}[thm]{Remark}

\newcommand\be{\begin{equation}}
\newcommand\ee{\end{equation}}
\newcommand\bea{\begin{eqnarray}}
\newcommand\eea{\end{eqnarray}}
\newcommand\bc{\begin{center}}
\newcommand\ec{\end{center}}
\newcommand\ba{\begin{array}}
\newcommand\ea{\end{array}}

\newcommand{\R}{\ensuremath{\mathbb{R}}}
\newcommand{\C}{\ensuremath{\mathbb{C}}}

\newcommand{\ga}{\alpha}                  

\newcommand{\gd}{\delta}

\newcommand{\gO}{\Omega}

\newcommand{\dgO}{\partial\Omega}

\newcommand{\gep}{\epsilon}

\newcommand{\esig}{e^{\sigma}}
\newcommand{\esigz}{e^{\sigma_0}}

\newcommand{\eg}{e^\gamma}
\newcommand{\egz}{e^{\gamma_0}}

\newcommand{\abs}[1]{\left| #1 \right|}
\newcommand{\norm}[1]{\left|\left|#1\right|\right|}

\begin{document}
\title[Stabilizing Inverse Problems]{Stabilizing Inverse Problems by Internal Data}
\author{Peter Kuchment$^1$ and Dustin Steinhauer$^2$}
\address{$^1$ Mathematics Department, Texas A\&M University, College Station, TX 77843-3368}
\address{$^2$ Inst. Appl. Math. and Comp. Sci., Texas A\&M University, College Station, TX 77843}
\ead{\mailto{kuchment@math.tamu.edu}, \mailto{dsteinha@math.tamu.edu}}

\begin{abstract} Several newly developing hybrid imaging methods (e.g., those combining electrical impedance or optical imaging with acoustics) enable one to obtain some auxiliary interior information (usually some combination of the electrical conductivity and the current) about the parameters of the tissues. This information, in turn, happens to stabilize the exponentially unstable and thus low resolution optical and electrical impedance tomography.

Various known instances of this effect have been studied individually. We show that there is a simple general technique (covering all known cases) that shows what kind of interior data stabilizes the reconstruction, and why. Namely, we show when the linearized problem becomes elliptic pseudo-differential one, and thus stable. Stability here is meant as the problem being Fredholm, so the local uniqueness is not shown and probably does not hold in such generality.
\end{abstract}

\maketitle

\section{Introduction}

The fast developing in the recent years \textbf{hybrid imaging methods} refer to a range of techniques in medical imaging in which different modalities are used in concert to benefit from the strengths of each while mitigating their individual weaknesses \cite{Ammari_book,bal2,Banff,Kuch_hybrid,KuKuEJM,KuKuElectr,WangPATbook,wang2}.  For example, ultrasound tomography provides high resolution, while not necessarily providing high contrast.  On the other hand, Electrical Impedance Tomography (EIT) and Optical Tomography (OT) can provide high contrast but are typically plagued by high instability and thus poor resolution \cite{BB1,bor,CIN}.  Acousto-Electric Tomography (AET), also called Ultrasound Modulated Electrical Impedance Tomography (UMEIT) \cite{Ammari_book,Ammari_AET,Kuch_hybrid,KuKuSynt,KuKuAET,ScherAET,wang}, uses focused (physically or synthetically) ultrasound waves to perturb the conductivity inside the object of interest, and these perturbations can be measured using the EIT techniques. A similar approach is used in UMOT (Ultrasound Modulated Optical Tomography), which combines optical tomography with acoustics \cite{AllmBang,BalSchAOT,Nam,NamDob,wang2}.
The MREIT and CDI/CDII techniques use various combinations of EIT and MRI \cite{Nach1,Nach2,Nach3,SeoSIAM,SeoMREIT}. What is common for these (physically rather different) techniques is that after some manipulations (see \cite{Ammari_book,bal2,bal1,bal3,bal4,BalRenQPAT,BalSchAOT,BalUhl,ScherAET,KuKuAET,Nach1,Nach2,Nach3,SeoSIAM,SeoMREIT} for details), the values through the interior of the object of a function (often of the form $\sigma(x)\abs{\nabla u(x)} ^p$) can be obtained.  Here $\sigma$ is the conductivity to be determined,
and $u$ is the corresponding electric potential.  The problem then becomes to determine $\sigma$ from this
interior data.

  Surveys of the recent results on hybrid methods can be found, for instance, in \cite{bal2,Kuch_hybrid,Ammari_book,KuKuEJM}.

It has been observed that a common feature of these methods is that they provide significantly better stability and resolution than more conventional EIT and OT techniques. The opinion has been expressed by several experts that some meta-statement should exist that claims that ``appropriate'' interior information stabilizes the exponentially unstable problems such as EIT or OT\footnote{The main difficulty in EIT and OT is that a signal by the time it reaches the boundary has largely ``forgotten'' where it originated from. Thus the idea is that having some quantity $F(\sigma(x), u(x),\nabla u(x))$ attached to an arbitrary interior location $x$ should compensate for this ``lack of memory.''}. Our aim is to provide a version (in fact, several versions of various generality) of such a statement. Doing this, we can address the stability of different hybrid methods with internal data in a unified way.

Our considerations are local, i.e. in a neighborhood of a known smooth background. We thus use linearization. If one could prove that the Fr\'{e}chet derivative of the corresponding non-linear mapping exists, has zero kernel, and is a (semi-) Fredholm operator in some scale of (say, Sobolev) Banach spaces, by standard functional analysis arguments this would have implied local uniqueness and stability of the non-linear problem. In this text, we provide very general statements on (semi-) Fredholm property of such Fr\'{e}chet derivatives (see, e.g., \cite{Krein,Kato,GoKr,ZKKP} for basics on Fredholm and semi-Fredholm operators). This is done by reducing them to elliptic pseudo-differential operators. This achieves  necessary stability estimates {\bf modulo the possible finite-dimensional kernel}. The authors doubt that the uniqueness claim (i.e., absence of non-trivial kernel) can be made in such wide generality in which we obtain ellipticity and thus Fredholm property. Thus, the absence of kernel should be dealt with on case-by-case basis (as it has been done before). However, existence of Fredholm property shows that the interior data does have stabilizing effect on the initially exponentially unstable problem.






Let us describe the structure of the article.
In Section \ref{S:GS}, we describe the framework for later sections, and we state basic definitions and lemmas.
In Section \ref{S:ICP} we investigate inverse conductivity problems with different types of additional interior data.
We first consider the data of the form $\sigma(x)|\nabla u(x)|^p$, which is known to arise in a variety of hybrid problems \cite{Ammari_book,bal2,bal1,bal3,bal4,BalRenQPAT,BalSchAOT,BalUhl,ScherAET,KuKuAET,Nach1,Nach2,Nach3,SeoSIAM,SeoMREIT,Kuch_hybrid}. For $p\in(0,1)$ we show that ellipticity, and thus stability arises with a single set of the interior data. For $p>1$, two such sets are needed. Finally, we look at a rather general type of interior data of the form $F(\sigma(x),u(x),\nabla u(x))$ and obtain sufficient conditions under which one gets the problem stabilized.
In Section \ref{S:qpat} we treat a problem coming from the so-called Quantitative Photo-Acoustic Tomography (QPAT) \cite{BalJoll,bal4,BalRenQPAT,Cox3:2009,CoxQPAT,Cox4:2009,Kuch_hybrid,KuKuAET}, where the equation becomes of the diffusion type (with an absorption term) rather than just a divergence type equation as in the previous Section.

The main results of the paper are contained in Theorems \ref{the1}, \ref{T:F3}, \ref{irene}, and \ref{T:qpat}.

Sections \ref{S:lemmas} and \ref{S:remarks} contain the proofs of some technical lemmas used and final remarks and conclusions correspondingly.

\section{Preliminaries}\label{S:GS}

Let $\gO$ be a bounded open region in $\R^n$ with smooth boundary, and let $\gO'$ be an open region compactly contained within $\gO$. We will also need an intermediate domain $\gO''$:
$$
\gO'\Subset \gO'' \Subset \gO.
$$
We will frequently need to use two cutoff functions, $\chi_1\in C_0^\infty(\Omega)$ that is equal to $1$ in a neighborhood of $\overline{\Omega''}$ and $\chi_2\in C_0^\infty(\Omega'')$ that is equal to $1$ in a neighborhood of $\overline{\Omega'}$.

In Section \ref{S:ICP} we will address various inverse conductivity problems,
in which the goal is to determine the unknown log-conductivity $\sigma$ in the elliptic problem
\footnote{It will be more convenient for us to work with log-conductivity, which we denote $\sigma$, rather than with the true conductivity $\exp \sigma$.}
\be 
L_\sigma u=-\nabla \cdot (\esig\nabla u)=0
\ee
in $\Omega$ from some boundary data. For instance, in the EIT (Electrical Impedance Tomography) the data might be the whole Dirichlet-to-Neumann map on $\partial\Omega$
\cite{bor,cal,CIN,Uhlm_asteris,IOut1,IOut2,UhlmCald}.

We thus will have to work frequently with the Dirichlet boundary
value problem
\be\label{eit}
\cases{
-\nabla \cdot (\esig\nabla u)=0\\
 u|_{\dgO} = f .
}
\ee
It will be important for us that the solution of (\ref{eit}) depends on log-conductivity $\sigma$. To emphasize this dependence,
we will sometimes write $u(\sigma)$ for the solution to (\ref{eit}).

Analogously to \cite{cap},
we define the affine space of admissible log-conductivities as

\be L^\infty_{\mathrm{ad}}(\gO)=\{\sigma\in L^\infty(\gO)\:|\:\sigma|_{\gO\backslash\gO'}=0\}   \;. \ee

Functions in $ L^\infty_{\mathrm{ad}}(\gO)$
can be considered as defined on $\R^n$ by extending them
by zero.
We assume that $\sigma\in L^\infty_{\mathrm{ad}}(\gO)$ and $f\in H^{1/2}(\dgO)$ in (\ref{eit}),
so $u\in H^1(\gO)$ (see. e.g. \cite{bor}).

In Section \ref{S:qpat} we will study the following more general problem:
\bea\label{E:qpat}
\cases{
L_{\sigma,\gamma}u:=-\nabla \cdot (e^\sigma\nabla u)+e^\gamma u=0 \\
u|_{\dgO}=f.
}
\eea
The coefficients $\sigma$ and $\gamma$ are the log-diffusion and log-attenuation coefficients, respectively.

\begin{lem}\label{frdiff}
The map from $(\sigma,\:\gamma)$ to $u$, defined by (\ref{E:qpat}), is Fr\'{e}chet differentiable as a mapping from $L^\infty(\gO)\times L^\infty(\gO)$ to $H^1(\gO)$
at any point $(\sigma_0,\:\gamma_0)\in C^\infty(\gO)\times C^\infty(\gO)$.
\end{lem}
This fact is well known, but we supply in Section \ref{S:lemmas} its proof.
Infinite differentiability of the coefficients $(\sigma_0,\:\gamma_0)$ is in fact an overkill assumption here, but we will use (and thus prove) the lemma only in the smooth case.

In many hybrid imaging methods (see, e.g. \cite{Ammari_book,Ammari_AET,bal2,ScherAET,KuKuAET,Nach1,Nach2,Nach3,SeoSIAM,SeoMREIT}), internal information of the form
\begin{equation}\label{E:Fgeneral}
F(\sigma(x),u(x),\nabla u(x)),
\end{equation}
where $u$ is the solution of (\ref{eit}) corresponding to some specific boundary data $f$, can be derived from the measured data. Thus, the next goal is to recover the log-conductivity function $\sigma(x)$ from the knowledge of $F(\sigma(x),u(x),\nabla u(x))$ for all $x\in\Omega$ (which explains the name ``interior data''). Since, for a fixed boundary Dirichlet data $f$, the functions $u$ and $\nabla u(x)$ are determined by $\sigma$, we can consider $F$ as a non-linear operator acting on $\sigma$: $F:\sigma \mapsto F(\sigma)$.


The pseudo-differential technique that we use makes the statements simple and their proofs rather transparent. We make use of basic facts about pseudodifferential operators on $\R^n$ (see, for instance, \cite{shubin} or \cite[Ch. 7]{taylor}). We will use the standard symbol classes $S^m(\R^n)$ comprising smooth functions $a(x,\xi)$ that satisfy for any multi-indices $\alpha,\beta$ and for sufficiently large $|\xi|$ the estimates

\be \left|\partial^\alpha_x\partial^\beta_\xi a(x,\xi)\right|
\leq C_{\ga\beta}(1+|\xi|)^{m-|\beta|} \nonumber\ee
with some constants $C_{\ga\beta}$.
(Such symbol classes are often denoted $S^m_{1,0}(\R^n)$, but we will omit the subscripts as we will not be considering more general symbol classes.)
The corresponding classes of pseudo-differential operators, which we denote by
$\mathrm{OP}S^m(\R^n)$, are given by

\be a(x,D)u=\frac{1}{(2\pi)^n}\int\!\!\int a(x,\xi)e^{i(x-y)\cdot\xi}u(y)\:dy\:d\xi\nonumber\ee
where $a(x,\xi)\in S^m(\R^n)$ (with the standard regularization of this expression, see e.g. \cite{taylor,shubin}).

If a symbol $a(x,\xi)\in S^m(\R^n)$ satisfies

\be a(x,\xi)=a_m(x,\xi)+r(x,\xi), \nonumber\ee
where $r(x,\xi)\in S^{m-1}(\R^n)$, we will call $a_m(x,\xi)$ the {\it principal symbol}  of $a(x,\xi)\in S^m(\R^n)$.
The principal symbol is determined modulo $S^{m-1}(\R^n)$.
If, for some $R>0$, the principal symbol $a_m(x,\xi)$ satisfies the estimate

\be |a_m(x,\xi)|\geq C|\xi|^{m}\mbox{ for }|\xi|\geq R, \nonumber\ee
then the symbol $a(x,\xi)$ is called {\it elliptic}, and the corresponding operator
$a(x,D)$ is called {\it elliptic} as well.

We will also need some facts about (square) matrix pseudo-differential operators. Let $A_{i,j}(x,\xi)$ for $i,j=1, ... , p$ be a matrix of classical symbols of pseudo-differential operators. Suppose that there exist two $p$-tuples $\{s_1, ... ,s_p\}$ and $\{t_1, ... ,t_p\}$ of real numbers such that $A_{i,j}\in S^{s_i+t_j}(\R^n)$. Let also $A_{i,j}^0\in S^{s_i+t_j}$ be their principal symbols. The system $A=\{A_{i,j}\}$ is said to be {\it elliptic in the Douglis-Nirenberg (DN) sense}\footnote{Sometimes it is also called Agmon-Douglis-Nirenberg ellipticity. An equivalent, although differently formulated, notion was also introduced by L.~Volevich (e.g., \cite[Ch. 9]{Wloka}).} (see \cite{ADN} or \cite[Ch. 7]{taylor}, \cite[Ch. 9]{Wloka}) if the determinant $\det \left(A^0_{i,j}(x,\xi)\right)$ does not vanish for $|\xi|>R$, for a suitable $R$.

\section{Stability in Inverse Conductivity Problems with Internal Data}\label{S:ICP}

Here, we address functionals $F$ of the more specific form
\begin{equation}\label{E:Fp}
F(\sigma)=\esig |\nabla u(\sigma)|^p.
\end{equation}
Several hybrid imaging methods provide this kind of internal data  (see, e.g., \cite{Ammari_book,bal2,bal1,bal3,bal4,BalRenQPAT,BalSchAOT,BalUhl,ScherAET,KuKuAET,Nach1,Nach2,Nach3,SeoSIAM,SeoMREIT}). For example, it arises with $p=2$ in AET (Acousto-Electric Tomography, also called sometimes UMEIT (Ultrasound Modulated Electrical Impedance Tomography), or Impediography).

When $0<p\leq2$ and $\sigma\in L^\infty_{\mathrm{ad}}(\gO)$, one concludes that $F(\sigma)(x)$ belongs to
$L^1(\gO)$.
Let us consider $F$ as a non-linear mapping from $L^\infty_{\mathrm{ad}}(\gO)$ to $L^1(\gO)$. As such, it is Fr\'{e}chet differentiable at any smooth log-conductivity $\sigma_0$,
as long as the corresponding solution $u(\sigma_0)$ has a gradient that is bounded below by a positive constant.
This is a direct consequence of the following lemma:

\begin{lem}\label{L:diff_general}
Let  $F(y,z,w)$ be a function of three variables that is smooth when $y, z\in\R$, and $w>0$.
Assume that $F$ satisfies the bound

\be \label{cofy}
\abs{F(y,z,w)}\leq C(y)(z^2+w^2)\;,
\ee
where $C(y)$ depends continuously on $y$.
Then the mapping $\sigma \Rightarrow F(\sigma, u, \abs{\nabla u})$
is Fr\'{e}chet differentiable at the smooth background log-conductivity $\sigma_0$, as a mapping from $L^\infty_\mathrm{ad}(\gO)$ to $L^1(\gO)$.
\end{lem}
The proof is provided in Section \ref{S:lemmas}.

Since the Fr\'{e}chet derivative exists, it can be found by a formal linearization.  Consider a small perturbation of $\sigma_0$:

\be \label{linearize}
\begin{array}{l}
\sigma=\sigma_0+\gep\rho \\
u=u_0+\gep v+o(\gep),
\end{array}
\ee
where $u_0=u(\sigma_0)$. A simple substitution, as in \cite{cap}, shows that $v\in H^1_0(\gO)$ solves the boundary value problem
\be \label{adj}
\cases{
L_{\sigma_0}v=\nabla\cdot(\rho\esigz\nabla u_0)\\
v|_{\dgO}=0 .
}
\ee
One notices that the dependence of $v$ on $\rho$ is linear. We will indicate it as $v(\rho)$.

The Fr\'{e}chet derivative $dF$ of $F$, which we will denote by $dF$, is a linear bounded operator from $L^{\infty}_{\mathrm{ad}}(\gO)$ to $L^1(\gO)$. Applying the chain rule and (\ref{adj}) to $\esig |\nabla u(\sigma)|^p$, one finds $dF$ as

\be \label{deeeff}
dF(\rho) = \rho\esigz|\nabla u_0|^p + p\esigz\frac{\nabla u_0\cdot\nabla v(\rho)}{|\nabla u_0|^{2-p}} \;.
\ee
We introduce a cutoff version of $dF$, which extends to a pseudo-differential operator on $\R^n$.
Let $\chi_1$ be a smooth cutoff function
supported in $\gO$ which is identically equal to 1 on a neighborhood of $\overline{\gO''}$.
We define an operator $A$ mapping $L^{\infty}(\gO)$ to $L^1(\R^n)$ by
\begin{equation}\label{E:A}
A(\rho)=\chi_1\:dF(\chi_1\rho).
\end{equation}

Because of the presence of the cutoff by $\chi_1$ before applying $dF$, the operator $A$ has a natural extension
to $L^\infty(\R^n)$.

In order to show that $A$ is a pseudodifferential operator, we analyze equation (\ref{adj}).
  The expression $\nabla\cdot(\chi_1\rho\esigz\nabla u_0)$, as a differential operator acting on $\rho$, has principal symbol $i\esigz\chi_1\xi\cdot\nabla u_0$.
This operator, when acting on functions that vanish outside of $\gO'$, does not depend on the choice of $\chi_1$.

The principal symbol of the elliptic differential operator $L_{\sigma_0}$ is $\esigz |\xi|^2$. Hence,
$L_{\sigma_0}$ has a pseudo-differential parametrix $P\in \mathrm{OP}S^{-2}(\R^n)$ with
principal symbol $(\esigz|\xi|^2)^{-1}$.  This means that $L_{\sigma_0}P=I+S$,
where $I$ is the identity operator on $\R^n$ and $S$ is a smoothing operator on $\R^n$, and analogously for $PL_{\sigma_0}$.
Let us define the following function:

\begin{equation*}
w:=P(\nabla\cdot(\chi_1\rho\esigz\nabla  u_0)).
\end{equation*}
Then we have
\bea
L_{\sigma_0}(v-w)&=&\nabla\cdot(\rho\esigz\nabla u_0)-L_{\sigma_0}P(\nabla\cdot(\chi_1\rho\esigz\nabla  u_0))\nonumber\\
&=& \nabla\cdot(\rho\esigz\nabla u_0)-(I+S)(\nabla\cdot(\chi_1\rho\esigz\nabla  u_0)) \nonumber\\
\label{vanillasmoothie} &=& \nabla\cdot\left( (1-\chi_1)\rho\esigz\nabla u_0\right)+
S(\nabla\cdot(\chi_1\rho\esigz\nabla  u_0)).
\eea
The expression in (\ref{vanillasmoothie}) is a smooth function, so $L_{\sigma_0}(v-w)\in C^\infty(\gO)$.
By elliptic regularity, $v\equiv w$ mod $C^\infty(\gO)$.  Because of this equivalence, the mapping
$\rho\mapsto v$ is a pseudo-differential operator modulo infinitely smoothing operators on $\gO$.  All other operations in equation (\ref{deeeff}) are simply multiplication operators,
so we see that after multiplying by $\chi_1$, $A$ is a pseudodifferential operator on $\R^n$.


Let $A_0(x,\xi)$ denote the principal symbol of $A$, so that $A(x,\xi)=A_0(x,\xi)+R_{-1}(x,\xi)$
where $R_{-1}$ is a symbol of order $-1$.  The symbol $A_0$ is easily derived.
The principal symbol of a composition of operators is the product of the
individual principal symbols.  Applying this to the composition of the operators

\bea \rho&\mapsto&\nabla\cdot(\chi_1\rho\esigz\nabla u_0) \nonumber\\
u&\mapsto& P(u) \nonumber\eea
shows that the mapping $\rho\mapsto v(\rho)$ has a principal symbol given by

\begin{equation*} \left(i\esigz\chi_1\xi\cdot\nabla u_0\right)(\esigz|\xi|^2)^{-1}=\frac{-i\chi_1 \xi\cdot\nabla u_0}{|\xi|^2} \;.\nonumber\end{equation*}
From equation (\ref{deeeff}) we then find that $A_0$ is given by

\be \label{E:A0}
A_0(x,\xi)
=\chi_1^2\esigz|\nabla u_0|^p + p\chi_1^2\esigz\frac{(i\xi\cdot\nabla u_0)^2}{|\nabla u_0|^{2-p}|\xi|^2}.
\ee

Let $\theta$ denote the angle between $\xi$ and $\nabla u_0$.
 Then the principal symbol $A_0$ becomes
 \begin{equation}\label{E:Aprinc}
 A_0(x,\xi)=\chi_1^2\esigz|\nabla u_0|^p(1-p\cos ^2 \theta)
 \end{equation}
 on $\gO'$.  Since we are interested in determining under what conditions $A$ is elliptic,
we are therefore motivated to consider the case when $p<1$ separately from $p\geq 1$.

\subsection{The $p<1$ case}\label{SS:p<1}

\begin{thm}\label{the1}If $p<1$, then
\begin{enumerate}
\item $A(x,D)$, as defined above, is a pseudo-differential operator of order zero, which is elliptic in
a neighborhood of $\overline{\gO''}$;
\item $dF$, as an operator acting in $L^2(\gO')$, is Fredholm;
\item Let $K$ be the kernel of $dF$ as an operator on $L^2(\gO')$, and let $R$ be its range ($K$ is finite-dimensional and $R$ is of a finite co-dimension). Then $dF$, considered as an operator from $L^2(\gO')/K$ onto $R$, is a topological isomorphism, i.e.
 there exists a constant $C$ such that for $\rho\in L^2(\gO')$,
\be \label{Freddie}
\frac{1}{C} \norm{dF( \rho)}_{L^2(\gO')} \leq \norm{\rho}_{L^2(\gO')/ K} \leq C \norm{dF(\rho)}_{L^2(\gO')} \;. \ee
\end{enumerate}
\end{thm}

\begin{rek}
\begin{enumerate}\indent
\item When we consider $dF$ as an operator from $L^2(\gO')$ to itself, we are really considering the operator $T\circ\:dF$, where $T:L^2(\R^n)\rightarrow L^2(\gO')$ is the restriction operator. The same goes for $A$.
\item Sometimes we need to consider $T$ as acting into $L^2(\R^n)$, in which case the action of $T$ is simply multiplication by the characteristic function of $\gO'$. The operator $T$ also maps the
corresponding Sobolev spaces:
\be
T:H^s(\R^n)\rightarrow H^s(\gO').
\ee
\end{enumerate}
\end{rek}

{\bf Proof.}
Since $A_0(x,\xi)=\chi_1^2\esigz|\nabla u_0|^p(1-p\cos ^2 \theta)$ and $p<1$, $A$ is elliptic on a neighborhood of $\overline{\gO''}$. This proves the first claim of the theorem.

Let $\chi_2$ be a smooth cutoff function supported in the $\Omega''$ that is equal to 1 on a neighborhood
of $\gO'$.  We define a symbol $Q(x,\xi)\in S^0(\R^n)$ by setting

\be Q(x,\xi)=\frac{\chi_2(x)}{A_0(x,\xi)} \;.\ee
Then, letting $I$ be the identity operator, $Q(x,D)A(x,D)-\chi_2(x) I$ is a pseudodifferential operator of order $-1$ on $\R^n$.

We claim that $TQ$ is a left inverse for $TA$ modulo compact operators on $L^2(\gO')$.  Indeed,
\be \label{ops}
TQTA-I=(TQA-I)+TQ(T-I)A .
\ee
Since the pseudo-differential operator $Q(x,D)A(x,D)-\chi_2(x) I$ is of order $-1$ on $\R^n$, the function $(TQA-I)f$ belongs to $H^1(\gO')$ for any $f\in L^2(\gO')$.
In addition, the function $(T-I)Af$ is equal to zero on $\gO'$, so by the microlocal
property of pseudodifferential operators, $TQ(T-I)Af$ is a smooth function on $\gO'$.  Thus the righthand side of equation (\ref{ops}) is a bounded operator from $L^2(\gO')$ into $H^1(\gO')$.  The imbedding operator of $H^1(\gO')$ into $L^2(\gO')$ is compact, proving the claim.

We also observe that
\be A(x,D)Q(x,D)\equiv Q(x,D)A(x,D)
\ee
modulo operators of order $-1$, so $TAQ-I$ maps $L^2(\gO')\rightarrow H^1(\gO')$.  Therefore $TA$ is also
right-invertible modulo compact operators.  This means that $TA$ is a Fredholm operator on $L^2(\gO')$.
Then $T\circ dF$ is a Fredholm operator on $L^2(\gO')$ as well, because $dF=A$ as operators on $L^2(\gO')$. This proves the second claim of the theorem.

The 3rd claim of the theorem is an immediate consequence of the 2nd one.
$\square$



\subsection{The $p\geq1$ case}\label{SS:p>1}
If $p\geq 1$, the formula (\ref{E:Aprinc}) shows that there are directions $\xi$ at each point $x$ in which the principal symbol $A_0(x,\xi)=\chi_1^2\esigz|\nabla u_0|^p(1-p\cos ^2 \theta)$ vanishes.
  In order to make the problem elliptic, we need to assume availability of two measurements. Namely,
  for two different boundary conditions $f_1,\: f_2$ in (\ref{eit}), let $u^{(i)}$ for $i=1\:,2$ be the corresponding solutions:
\be
\cases{\label{burpie}
L_\sigma u^{(i)}=0\\
	 u^{(i)}|_{\dgO}= f_i.
}
\ee
Assume that we are given knowledge, for all $x\in\gO$, of the following sets of internal data:
\be \label{E:2data}
\cases{
F_{11}=\esigz\abs{\nabla u^{(1)}}^p\\
F_{22}=\esigz\abs{\nabla u^{(2)}}^p\\
F_{12}=\esigz\left| \nabla u^{(1)}\cdot \nabla u^{(2)}\right|^{p/2}.
}
\ee
Such functionals have been extracted from the measured data in hybrid imaging methods, see
for example \cite{bal3,cap,KuKuAET}.

Our (vector) internal measurement function will be now
\bea
F(\sigma):=\left(
     \begin{array}{c}
       F_{11}(\sigma) \\
       F_{22}(\sigma) \\
       F_{12}(\sigma) \\
     \end{array}
   \right).
\eea

  We again consider small perturbations of a smooth background log-conductivity $\sigma_0$ as in (\ref{linearize}),
and let $u^{(i)}_0$ be the corresponding solutions of (\ref{burpie}) with $\sigma_0$ as the log-conductivity:
\be
\cases{\label{burpie_0}
L_{\sigma_0} u_0^{(i)}=0\\
	 u^{(i)}_0|_{\dgO}= f_i.
}
\ee
We need to assume in addition that the gradients of $u^{(1)}_0$ and $u^{(2)}_0$ are nowhere parallel:
\be \label{notparallel}
u^{(1)}_0\not{\parallel} u^{(2)}_0.
\ee
This is known to be possible \cite{ale} in $2D$ under an appropriate choice of $f_1$ and $f_2$.
However, as shown in \cite{lau}, in $3D$, it is not always possible to choose boundary conditions such that
(\ref{notparallel}) is satisfied. Condition (\ref{notparallel}) means in particular that the gradients are nonvanishing.

\begin{lem}\label{L:diffF3}
\indent
\begin{enumerate}
\item  As in the previous sub-section, functionals $F_{11}$ and $F_{22}$ are Fr\'{e}chet differentiable with respect to $\sigma$ as mappings from $L^\infty_{ad}$ to $L^1$.
\item When $p<2$, the mapping $F_{12}$ is Fr\'{e}chet differentiable if $u^{(1)}_0$ and $u^{(2)}_0$
satisfy
\be \label{gogo}
\nabla u^{(1)}_0 \cdot \nabla u^{(2)}_0 \geq \ga
\ee
for some $\ga>0$.  If this condition fails at some points, it is still true that
the functional $\phi F_{12}$ is Fr\'{e}chet differentiable if $\phi$ is a smooth cutoff function
and (\ref{gogo}) is satisfied on its support.
\item When $p=2$, mapping $F_{12}$ is Fr\'{e}chet differentiable without condition (\ref{gogo}).
\end{enumerate}
\end{lem}

\begin{rek}\indent
\begin{itemize}
\item We will handle the case $1<p<2$ and note that the case $p=2$ follows by the same argument, with the part concerning the smoothness of $F_{12}$ omitted.
\item For the case $p=1$, the functional $F_{12}$ is simply not needed, though the following arguments still apply.
Thus the full range $1\leq p\leq2$ will be covered.
\end{itemize}
\end{rek}

The proof of this lemma can be found in Section \ref{S:lemmas}.

Given a vector $\xi$, let $\theta_1$ be the angle between $\xi$ and $\nabla u^{(1)}_0$,
$\theta_2$ be the angle between $\xi$ and $\nabla u^{(2)}_0$,
 and $\theta$ be the angle between $\nabla u^{(1)}_0$ and $\nabla u^{(2)}_0$.
As in the case when $p<1$, we define cutoff versions of $dF_{ij}$ by setting $A_{ij}=\chi_1 dF_{ij}\chi_1$.
The principal symbols of $A_{11}$ and  $A_{22}$
near $\overline{\gO'}$,  calculated in the same manner as before, are

\bea
A_{11}(x,\xi)&=& \chi_1^2\esigz|\nabla u^{(1)}_0|^p(1-p\cos ^2 \theta_1), \nonumber\\
A_{22}(x,\xi)&=& \chi_1^2\esigz|\nabla u^{(2)}_0|^p(1-p\cos ^2 \theta_2). \nonumber
\eea
The principal symbol of $A_{12}$ at points where $\nabla u^{(1)}_0$ and $\nabla u^{(2)}_0$
are not perpendicular is also easily derived from the formula

\bea dF_{12}&=&\esigz\rho\abs{\nabla u^{(1)}_0\cdot\nabla u^{(2)}_0}^{p/2} \nonumber\\
&+&\esigz\frac{p}{2}\abs{\nabla u^{(1)}_0\cdot\nabla u^{(2)}_0}^{p/2-1}
\left(\nabla u^{(1)}_0\cdot\nabla v^{(2)}(\rho)+\nabla u^{(2)}_0\cdot\nabla v^{(1)}(\rho)\right) \eea
with

\be \label{adj}
\cases{
L_{\sigma_0}v^{(i)}=\nabla\cdot(\rho\esigz\nabla u_0)\\
v^{(i)}|_{\dgO}=0 ,
}
\ee
which readily follows from formal linearization.  Hence,

\be \label{f12}
A_{12}(x,\xi)= \chi_1^2\esigz\abs{\nabla u^{(1)}_0}^{p/2} \abs{\nabla u^{(2)}_0}^{p/2}|\cos^{p/2}\theta|
\left(1-\frac{p\cos\theta_1\cos\theta_2}{|\cos\theta|}\right) \;.\ee

If, near each point $(x,\xi)\in\overline{\gO'}\times\R^n\backslash 0$,
at least one of these symbols is non-vanishing,
then the (vector) operator $(A_{11}\;A_{22}\;A_{12})^t$ is overdetermined elliptic.

  We notice that $A_{ii}(x,\xi)$ vanishes when $\cos\theta_i=\pm\frac{1}{\sqrt p}$.
  Near points $(x,\xi)$ where both $A_{11}(x,\xi)$ and $A_{22}(x,\xi)$ vanish, the symbol $A_{12}(x,\xi)$ will have to save the situation. Near such a point,  i.e. where $|\cos\theta_1|=|\cos\theta_2|=\frac{1}{\sqrt p}$, a simple trigonometric estimate shows that the expression $|\cos\theta|$ is separated from zero, and thus the symbol $A_{12}(x,\xi)$ is smooth. The non-vanishing of $A_{12}(x,\xi)$ then boils down to $\cos\theta\neq 1$ at the point $x$, which is guaranteed by (\ref{notparallel}).

Let us define three symbols homogeneous of order 0 on $\Omega\times\R^n$,
 $\psi_{11}(x,\xi)$, $\psi_{22}(x,\xi)$, and $\psi_{12}(x,\xi)$, such that

\be \Psi(x,\xi)= \psi_{11}(x,\xi)A_{11}(x,\xi)+\psi_{22}(x,\xi)A_{22}(x,\xi)+\psi_{12}(x,\xi)A_{12}(x,\xi) \ee
is a non-vanishing symbol of order 0 near $\overline{\gO'}\times\R^n\backslash 0$, and therefore bounded away from zero
by the compactness of the cosphere bundle of $\overline{\gO'}$.
The above arguments imply that $\psi_{12}$
can be taken to be zero near points where condition (\ref{gogo}) is not satisfied.
(For $p=1$ we do not need the functional $F_{12}$, so we set $\psi_{12}\equiv0$ in this case.)
The operators $\psi_{ij}(x,D)$ with symbols $\psi_{ij}(x,\xi)$
are bounded on $L^2(\Omega)$.  This means that the operator norm of $\Psi(x,D)$ can be controlled by the sum of the
operator norms of $A_{ij}(x,D)$.
We thus have

\begin{thm}\label{T:F3}

\begin{enumerate}
\indent
\item  For $1\leq p\leq 2$, the operator
$$
dF: L^2(\Omega') \Rightarrow \{L^2(\Omega')\}^3
$$
is semi-Fredholm with a possible finite dimensional kernel $K$.
\item Letting $K$ be the finite-dimensional kernel of $dF$, the estimate

\be \label{E:esimate}
 C^{-1}\norm{\rho}_{L^2(\gO')/ K}
\leq \norm{dF( \rho)}_{\{L^2(\gO')\}^3}
\leq C \norm{\rho}_{L^2(\gO')/ K}
\ee
holds for $\rho\in L^2(\gO')$ and some constant $C$.
\end{enumerate}
\end{thm}

\begin{rek} As we have already mentioned, such functionals arise naturally and have been studied previously.  When $p=2$, a similar local stability
estimate was proved in \cite{KuKuAET, cap} in the space $C^{1,\ga}(\gO')$,
and in \cite{bal3} a global estimate was established in $W^{1,\infty}(\gO)$.  A stability estimate for a single functional of the form $e^\sigma|\nabla u|^2$
was also established in \cite{bal1} on a part of $\gO$.  In the case $p=1$, inversion procedures and reconstructions for a single functional
were obtained in \cite{Nach1,Nach2,Nach3}.
\end{rek}

{\bf Proof.}
 As in the proof of Theorem \ref{the1}, 
let

\be Q(x,\xi)=\frac{\chi_2(x)}{\Psi(x,\xi)} \;.\nonumber\ee
Then the operator $Q(x,D)$ lies in $\mathrm{OP}S^0(\R^n)$. Arguments identical
to the ones in the proof of Theorem \ref{the1} show that
$T\Psi(\Psi_{11},\Psi_{2},\Psi_{12})$ is left regularizer\footnote{Operator $B$ is a \textbf{left regularizer} to operator $A$, if $BA-I$ is a compact operator \cite{ZKKP,Krein}.} for the operator
\begin{equation}\label{E:left}
\left(
  \begin{array}{c}
    F_{11} \\
    F_{22} \\
    F_{12} \\
  \end{array}
\right): L^2(\gO')\mapsto (L^2(\gO'))^n.
\end{equation}
Thus, the operator in (\ref{E:left}) is semi-Fredholm with a finite dimensional kernel (and infinite-dimensional co-kernel), which implies all the statements of the theorem.
$\square$

\subsection{More general interior data}\label{SS:general}

We next consider a single, rather general functional and formulate a sufficient condition for
the corresponding linearized problem being elliptic (and thus Fredholm).

Let  $F(y,z,w)$ be a function of three variables satisfying the conditions of Lemma \ref{L:diff_general}; that is, $F$ smooth when $y, z\in\R$, and $w>0$,
and satisfies the bound

$$
\abs{F(y,z,w)}\leq C(y)(z^2+w^2)\;,
$$
where $C(y)$ depends continuously on $y$.


The Fr\'{e}chet derivative can be derived by a formal calculation as before, and is given by

\be \label{asdf} dF(\rho)=\frac{\partial F}{\partial y}\rho + \frac{\partial F}{\partial z}v
+\frac{\partial F}{\partial w}\frac{\nabla u_0\cdot\nabla v}{\abs{\nabla u_0}}\;. \ee

As calculated before, the principal symbol of
the operator mapping $\rho$ to $v$ is
$i\xi\cdot\nabla u_0|\xi|^{-2}$, which is of order $-1$.  The middle term on the right hand side
in expression (\ref{asdf}) is therefore of lower order than the other two
terms on the right hand side and does not influence the ellipticity of the overall principal symbol.
Also as before, we set $A=\chi_1 dF$.
From the other two terms in (\ref{asdf}) we find that the principal symbol of $A$ is

\be A_0(x,\xi)= \frac{\partial F}{\partial y} - \frac{\partial F}{\partial w}
\frac{(\xi\cdot\nabla u_0)^2}{|\nabla u_0||\xi|^2} =
\frac{\partial F}{\partial y} - \frac{\partial F}{\partial w}|\nabla u_0|\cos ^2\theta \ee
near $\overline{\gO'}$ (we omit the cutoff function $\chi_1$ because it is identically equal to 1 there).
This leads to a sufficient condition for the ellipticity of $A_0$:

\begin{thm}\label{irene}
If
\be \label{lkjh}
\left| \frac{\partial F(\sigma_0,u_0,|\nabla u_0|)}{\partial y}\right|
> \norm{\nabla u_0}_{L^\infty(\gO)}\left|\frac{\partial F(\sigma_0,u_0,|\nabla u_0|)}{\partial w}\right| \ee
pointwise in a neighborhood of $\overline{\gO}$, then
\begin{enumerate}
 \item $A(x,\xi)$ is elliptic of order 0 on a neighborhood of $\overline{\gO'}$;
 \item $dF$ as a Fredholm operator in $L^2(\gO')$;
 \item Letting $K$ be the finite-dimensional kernel of $ dF$, the estimate

\be \label{guestimate}
\frac{1}{C} \norm{dF( \rho)}_{L^2(\gO')} \leq \norm{\rho}_{L^2(\gO')/ K} \leq C \norm{dF(\rho)}_{L^2(\gO')} \ee
holds for $\rho\in L^2(\gO')$.
\end{enumerate}
\end{thm}
\begin{rek}
\begin{enumerate}
 \item The functional $F(\sigma)=\sigma|\nabla u|^p$ satisfies assumption (\ref{lkjh}) for $p<1$, so Theorem \ref{irene}
 generalizes Theorem \ref{the1}.
 \item If $\sigma_0\equiv c$ for some constant $c$, we can have $|\nabla u_0|\equiv 1$ by selecting an appropriate boundary
condition, e.g. $f=x_1$.  Then (\ref{lkjh}) just means
$\left| \frac{\partial F}{\partial y}\right| > \left|\frac{\partial F}{\partial w}\right|$ pointwise on $\gO$.
\end{enumerate}
\end{rek}

{\bf Proof.}
We may assume that $\chi_1\equiv1$ on the region where condition (\ref{lkjh}) is satisfied.
If (\ref{lkjh}) holds, then there exists $\gd>0$ such that

\be \left| \frac{\partial F}{\partial y}\right| - \left|\frac{\partial F}{\partial w}\right| \norm{\nabla u_0}_{L^\infty(\gO)} >\gd>0\nonumber\ee
on $\overline{\gO}$ by the continuity of the terms involved.
This implies that $A_0(x,\xi)>\gd$, and so $A$ is elliptic on the same neighborhood of $\gO'$.
By setting

\be Q(x,\xi)=\frac{\chi_2(x)}{A_0(x,\xi)}\nonumber\ee
and composing with the restriction operator $T$,
we obtain an operator $T\circ Q(x,D)$ which is a left and right inverse for $TA$ on $L^2(\gO')$
(and hence for $T\circ dF$ as well) modulo compact operators.  Estimate (\ref{guestimate}) then follows.$\square$


\section{Stability in Quantitative Photo-Acoustic Tomography}\label{S:qpat}

The standard model for diffusive regime photon propagation in biological tissues is

\bea\label{qpat}  \cases{
L_{\sigma,\gamma}u:=-\nabla \cdot (e^\sigma\nabla u)+e^\gamma u=0 \\
u|_{\dgO}=f
}
\eea
(see, e.g. \cite{wang2}).
The coefficients $\sigma$ and $\gamma$ are the log-diffusion and log-attenuation coefficients, respectively. We will assume in this section that $\sigma \in L^\infty_{ad} (\Omega)$ and $\gamma \in H^1(\Omega)\cap L^\infty_{ad} (\Omega)$.

The PAT (Photo Acoustic Tomography) procedure done first, provides one with the values inside $\gO$ of the function

\be\label{mq}
F(x)=\Gamma (x) e^\gamma(x) u(x) .
\ee
Here, $\Gamma(x)$ is the so-called Gr\"{u}neisen coefficient\footnote{The Gr\"{u}neisen coefficient is in principle also not known, so one might want to include it as an unknown in the reconstruction procedure, e.g. \cite{BalRenQPAT}. We are not doing this here.}
describing the transfer of
electromagnetic into acoustic energy, which we assume here to be identically equal to 1.

This function is the initial data for the Quantitative Photo-Acoustic Tomography (QPAT), which strives to reconstruct the coefficients $\sigma$ and $\gamma$ from the data (\ref{mq}).

We will denote by $F_j(x)$, $j=1,2,\ldots ,J$ the internal data (\ref{mq}) that correspond to
solutions of (\ref{qpat}) with different boundary conditions $f_j$.

For a pair of such measurements $(F_1,\: F_2)$, according to Lemma \ref{frdiff}, the mapping
$(\sigma,\gamma)\Rightarrow F:=(F_1,\: F_2)$ is Fr\'{e}chet differentiable at smooth background coefficients $(\sigma_0,\gamma_0)$. The derivative can be computed formally as before:

\bea \sigma&=&\sigma_0+\gep\rho \nonumber\\
\gamma&=&\gamma_0+\gep\nu \nonumber\\
u^{(j)}&=&u_0^{(j)}+\gep v^{(j)}+o(\gep) \eea
where $\nu,\rho\in L^\infty_0(\gO)$. Substitution into (\ref{qpat}) shows that $v^{(j)}\in H^1(\gO)$ solves the boundary value problem

\be
\cases{
-\nabla\cdot(\esigz\nabla v^{(j)})+\egz v^{(j)}
= \nabla\cdot(\rho\esigz\nabla u_0^{(j)})-\nu \egz u_0^{(j)}\\
 v^{(j)}|_{\dgO}=0.
 }
 \ee
We thus find that the differential of the mapping $F_j$ is
\be
dF_j (\rho,\:\nu)=\nu u_0^{(j)} -L^{-1}_{\sigma_0,\gamma_0}\left(\nu u_0^{(j)}\right)+L^{-1}_{\sigma_0,\gamma_0}\left(\nabla\cdot(\rho\esigz\nabla u_0^{(j)})\right).
\label{E:dF}
\ee

A calculation similar to the one in the previous section shows that the operator $L_{\sigma_0,\gamma_0}$
has a parametrix with principal symbol

\be \frac{1}{\esigz|\xi|^2+\egz}\;, \nonumber\ee
which is equivalent to $(\esigz|\xi|^2)^{-1}$ modulo lower order terms.
Hence, according to (\ref{E:dF}), the matrix of the principal symbols of the operator $(\rho,\nu) \mapsto
 \chi_1 dF \chi_1$ is given by

\be A(x,\xi):= \chi_1
\left[ {\begin{array}{cc}
 \frac{i\xi\cdot\nabla u_0^{(1)}}{|\xi|^2} & u_0^{(1)}  \\
 \frac{i\xi\cdot\nabla u_0^{(2)}}{|\xi|^2} & u_0^{(2)}  \\
 \end{array} } \right]
\chi_1
\;,\ee
modulo lower order terms.

We consider here the Douglis-Nirenberg parameters
$$
s=(s_1,s_2)= (-1,-1), t=(t_1,t_2)= (0,1)
$$
and attempt to check the DN ellipticity of $dF$.

The determinant $\det A(x,\xi)$ is non-vanishing near $\overline{\gO'}$, if

\be \xi\cdot (u_0^{(1)}\nabla u_0^{(2)}-u_0^{(2)}\nabla u_0^{(1)}) \neq 0\;.\ee

Thus, ellipticity fails at the vectors $\xi$ orthogonal to the field $(u_0^{(1)}\nabla u_0^{(2)}-u_0^{(2)}\nabla u_0^{(1)})$. The natural idea is to have more measurements that would provide a basis of vector fields and thus preserve ellipticity. This leads to the following result:

\begin{thm}\label{T:qpat} Suppose that one has access to $2n$ measurements \
$$
(F_{1,1},F_{1,2}), \dots, (F_{n,1},F_{n,2}),
$$
such that the vector fields
$$
V_k:=u_0^{(k,1)}\nabla u_0^{(k,2)}
-u_0^{(k,2)}\nabla u_0^{(k,1)}
$$
for $k=1,\dots,n$ form a basis at each point $x$ in $\gO$. We define the operator $F$ as follows:
$$
F:=\left(F_{1,1},F_{1,2}, \dots, F_{n,1},F_{n,2}\right).
$$

Then the operator
$$
dF: L^2(\Omega')\bigoplus H^1_0(\Omega')\Rightarrow \{H^1(\Omega')\}^{2n}
$$
is semi-Fredholm with a finite dimensional kernel.

Letting $K$ be the (finite-dimensional) kernel of $dF$, the estimate

\bea 
\label{E:qpatstab}
 \frac{1}{C} \norm{\left(
                                  \begin{array}{c}
                                    \rho \\
                                    \nu \\
                                  \end{array}
                                \right)
 }_{H/K} \leq \norm{dF\left(
                                  \begin{array}{c}
                                    \rho \\
                                    \nu \\
                                  \end{array}
                                \right)}_{H^1(\gO'')^{2n}}\leq C \norm{\left(
                                  \begin{array}{c}
                                    \rho \\
                                    \nu \\
                                  \end{array}
                                \right)}_{H/K}
\eea
holds for some constant $C>1$.

Here we used the shorthand notation
$$
H:=L^2(\Omega')\bigoplus H^1_0(\Omega').
$$
\end{thm}

\begin{rek}
The assumptions we made on the vector fields $V_k$, as the reader could see, arose naturally. It is interesting to notice that they are the same that were also arising in the study of QPAT in \cite{BalUhl}, in which the CGO (complex geometric optics) solutions technique was used. The authors of \cite{BalUhl} derive a global estimate that is somewhat similar to (\ref{E:qpatstab}):
\be \label{E:Bal}
\norm{\delta\sigma}_{C^l(\gO)}+\norm{\delta\gamma}_{C^l(\gO)}
 \leq C
\mathop{\sum_{k=1,\dots,n}}_{j=1,2}
\norm{F_{k,j}(\sigma_1,\gamma_1)-F_{k,j}(\sigma_2,\gamma_2)}_{{C^{l+1}(\gO)}^{2n}}
\ee
for $l\geq2$, where $\delta\sigma=\sigma_1-\sigma_2, \delta\gamma=\gamma_1-\gamma_2$.

Under an additional convexity assumption on $\dgO$, the authors of \cite{BalUhl} also derive such an estimate with only two vector fields, when $l\geq 3$.

 The reader might also notice that in (\ref{E:qpatstab}), in comparison with (\ref{E:Bal}), different orders of smoothness are used, which allows us to get the two-sided estimate.
\end{rek}

{\bf Proof.}
The $2n\times 2$ matrix operator $\chi_2 dF\chi_2$ has the principal symbol
$$
A(x,\xi):= \chi_2\left(
  \begin{array}{cc}
    \frac{i\xi\cdot\nabla u_0^{(1,1)}}{|\xi|^2} & u_0^{(1,1)}  \\
 \frac{i\xi\cdot\nabla u_0^{(1,2)}}{|\xi|^2} & u_0^{(1,2)}  \\
 \dots & \dots\\
    \frac{i\xi\cdot\nabla u_0^{(n,1)}}{|\xi|^2} & u_0^{(n,1)}  \\
 \frac{i\xi\cdot\nabla u_0^{(n,2)}}{|\xi|^2} & u_0^{(n,2)}
  \end{array}\right)\chi_2.
$$
Here the principal symbol is understood in the Douglis-Nirenberg sense with parameters $s=(-1,-1,\dots,-1,-1), t=(0,1)$, and $\chi_2$, as before, is a smooth cutoff function that is equal to $1$ in a neighborhood of $\overline{\Omega'}$ and vanishes outside $\Omega''$.

Under the assumptions of the theorem,
at every $(x,\xi)\in \overline{\gO'}\times\left(\R^n\backslash 0\right)$
the symbol $A(x,\xi)$ is injective (since at least one of the square $2\times 2$ blocks
$$
A_k:=\left(
  \begin{array}{cc}
    \frac{i\xi\cdot\nabla u_0^{(k,1)}}{|\xi|^2} & u_0^{(k,1)}  \\
 \frac{i\xi\cdot\nabla u_0^{(k,2)}}{|\xi|^2} & u_0^{(k,2)}
  \end{array}\right)
$$
for $k=1,\dots,n$ is invertible).
Thus, the operator is overdetermined elliptic in DN sense.

Thus, there exists a left parametrix with the principal $2\times 2n$ symbol $B(x,\xi)$ with DN parameters $s=(0,-1), t=(1,1,\dots,1,1)$ (this is a well-known construction, which we indicate at the end of Section \ref{S:lemmas}). In other words, for the corresponding pseudo-differential operators $dF$ and $B$ one has $B\chi_2 dF\chi_2=I+K$ in a neighborhood of $\overline{\Omega'}$, where $K$ is a smoothing operator.

Let us agree to extend both $\rho\in L^2(\Omega')$ and $\nu\in H^1_0(\Omega')$ as equal to zero outside $\overline{\Omega'}$, without changing notations for these extended functions. Then, according to the DN parameters, $\chi_2 dF\chi_2\left(
     \begin{array}{c}
       \rho \\
       \nu \\
     \end{array}
   \right)
$ belongs to $(H^1_{0}(\Omega))^{2n}$, and the corresponding mapping from $L^2(\Omega')\bigoplus H^1_0(\Omega')$ to $(H^1_{0}(\Omega))^{2n}$ is continuous. This proves the right hand side inequality in (\ref{E:qpatstab}).

Similarly, $B$ acts continuously from  $(H^1_{comp}(\R^n))^{2n}$ to $L^2_{comp}(\R^n)\bigoplus H^1_{comp}(\R^n)$.

Then for the restriction of $\left(
     \begin{array}{c}
       \rho \\
       \nu \\
     \end{array}
   \right)$
   to $\Omega'$ we can write
  $$
  \begin{array}{l}
  \left(
     \begin{array}{c}
       \rho \\
       \nu \\
     \end{array}
   \right)=\chi_2(I+K)\left(
     \begin{array}{c}
       \rho \\
       \nu \\
     \end{array}\right) -\chi_2 K\left(
     \begin{array}{c}
       \rho \\
       \nu \\
     \end{array}
   \right)\\
   =\chi_2 B \chi_2 dF\chi_2\left(
     \begin{array}{c}
       \rho \\
       \nu \\
     \end{array}
   \right)-\chi_2K\left(
     \begin{array}{c}
       \rho \\
       \nu \\
     \end{array}
   \right)
   \end{array},
  $$
or
$$
  \begin{array}{l}
 (I+\chi_2 K) \left(
     \begin{array}{c}
       \rho \\
       \nu \\
     \end{array}
   \right)
   =\chi_2 B dF\left(
     \begin{array}{c}
       \rho \\
       \nu \\
     \end{array}
   \right)
   \end{array}.
  $$
  This gives us the estimate from above of the following kind:
  $$
  \|(I+\chi_2 K) \left(
     \begin{array}{c}
       \rho \\
       \nu \\
     \end{array}
   \right)\|_{L^2(\Omega ')\bigoplus H^1_0 (\Omega ')}\leq C \|dF\left(\begin{array}{c}
       \rho \\
       \nu \\
     \end{array}
   \right)\|_{H^1(\Omega'')^{2n}}.
  $$
  Since the operator $\chi_2 K$ is compact in $H=L^2(\Omega '')\bigoplus H^1_0 (\Omega '')$, the operator $I+\chi_2 K$ is Fredholm in $H$.

This implies the remaining statements of the theorem.
$\square$

\section{Proofs of some lemmas}\label{S:lemmas}

After proving Lemma \ref{frdiff}, we then proceed to give a proof of Lemma \ref{L:diff_general}.
Finally we prove Lemma \ref{L:diffF3}, making use of Lemma \ref{L:diff_general}.

\subsection{Proof of Lemma \ref{frdiff}}\label{SS:proof}
{\bf Proof.}
First of all, we reduce (\ref{E:qpat}) to a problem with homogeneous boundary conditions.
Let $E$ the operator of harmonic extension from $\partial\Omega$ to $\Omega$. Then, replacing the solution $u$ with $v_{\sigma,\gamma}:=u-Ef$, we reduce (\ref{E:qpat}) to
\begin{equation}\label{E:qpat_homog}
 \cases{
L_{\sigma,\gamma}v_{\sigma,\gamma}=f_{\sigma,\gamma}, \\
v_{\sigma,\gamma}|_{\dgO}=0,
}
\end{equation}
where
\begin{equation}
f_{\sigma,\gamma}:=\nabla \cdot (\esigz\nabla Ef)-\egz Ef\in H^{-1}(\Omega).
\end{equation}
Now the map $(\sigma,\gamma)\mapsto u_{\sigma,\gamma}$ factors as the composition of the following maps:
\begin{equation}\label{E:chain}
(\sigma,\gamma)\mapsto \{(\esig,\eg),f_{\sigma,\gamma}\}\mapsto L_{\sigma,\gamma}\mapsto L_{\sigma,\gamma}^{-1}\mapsto L_{\sigma,\gamma}^{-1}f_{\sigma,\gamma}+Ef.
\end{equation}

The first map in (\ref{E:chain}), is clearly Fr\'{e}chet differentiable as a mapping from
$L^\infty(\gO)\times L^\infty(\gO)$ to $L^\infty(\gO)\times L^\infty(\gO)\times H^{-1}(\Omega)$. Indeed, $f_{\sigma,\gamma}$ depends linearly and continuously on $(\sigma,\gamma)$.

The second map $(\esig,\eg)\mapsto L_{\sigma,\gamma}$ is linear and continuous from $L^\infty(\gO)\times L^\infty(\gO)$ to $L(H^1_0(\Omega),H^{-1}(\Omega))$.

The operator $L_{\sigma_0,\gamma_0}\in L(H^1_0(\Omega),H^{-1}(\Omega))$ is invertible (see the simplest case of this statement in \cite[Ch. 5, Proposition 1.1]{taylor} and general results in \cite{BerKrRo,Lions}); $L_{\sigma,\gamma}$ is thus invertible in a neighborhood of $(\gamma_0,\sigma_0)$ in $L^\infty(\gO)\times L^\infty(\gO)$.

The mapping
\bea L(H^1_0(\gO),H^{-1}(\gO))&\rightarrow&L(H^{-1}(\gO),H^1_0(\gO))
\eea
of taking inverse operator is known (e.g., \cite{ZKKP}) to be analytic on the domain of invertible operators. This implies differentiability of the last two mappings in (\ref{E:chain}) and thus proves the lemma.
$\square$

\subsection{Proof of Lemma \ref{L:diff_general}}

{\bf Proof.}
The function $F(\sigma(x), u(x), \abs{\nabla u(x)})$ lies in $L^1(\gO)$.
Indeed, we have $\abs{F(\sigma,u,|\nabla u|)}\leq C(\norm{\sigma}_{L^\infty(\gO)})(u^2+|\nabla u|^2)$, and $u$ and $|\nabla u|$ are both
square-integrable functions.

As a result of Lemma \ref{frdiff}, the map

\be \sigma\mapsto (\sigma,u,\nabla u) \ee
is Fr\'{e}chet differentiable from $L^\infty_\mathrm{ad}(\gO)\rightarrow L^\infty_\mathrm{ad}(\gO)\times L^2(\gO)\times \{L^2(\gO)\}^n$.
(The middle space could be taken to be $H^1(\gO)$ as before, but we will not need this here.)
We claim that the map

\bea L^\infty_\mathrm{ad}(\gO)\times L^2(\gO)\times \{L^2(\gO)\}^n\rightarrow L^1(\gO) \nonumber\\
(f,g,{\bf h})\mapsto F(f,g,|{\bf h}|) \nonumber\eea
is Fr\'{e}chet differentiable at $(f_0,g_0,{\bf h_0})\in L^\infty_\mathrm{ad}(\gO)\times L^2(\gO)\times \{L^2(\gO)\}^n$
when $f_0$, $g_0$, and ${\bf h_0}$ are bounded and smooth with $|{\bf h_0}(x)|\geq m$ for $x\in\gO$
and some positive constant $m$.
The boundedness of $f_0$, $g_0$, and ${\bf h_0}$ implies that for any multiindex $\ga$, the function

\be \label{partials} D^\ga F(y,z,w)|_{y=f_0(x),\:z=g_0(x),\:w=|{\bf h_0}(x)|} \ee
is a bounded function on $\gO$.

Let $(f,g,{\bf h})$ 
be a triple of functions in $L^\infty_0(\gO)\times L^2(\gO)\times \{L^2(\gO)\}^n$.  Consider the function
$E=E(x)$ defined on $\gO$ by

\bea E&=&F(f_0+f,g_0+g,|{\bf h_0}+{\bf h}|)-F(f_0,g_0,|{\bf h_0}|) \nonumber\\
&-&\nabla F(f_0,g_0,|{\bf h_0}|)
\cdot (f, g, \frac{{\bf h_0} \cdot{\bf h}}{{\bf |h_0|}}) \nonumber\\
&=&F(f,g,|{\bf h}|)-F(f_0,g_0,|{\bf h_0}|)-\frac{\partial F}{\partial y}(f_0,g_0,|{\bf h_0}|)f\nonumber\\
&-&\frac{\partial F}{\partial z}(f_0,g_0,|{\bf h_0}|)g-\frac{\partial F}{\partial w}(f_0,g_0,|{\bf h_0}|)\frac{{\bf h_0} \cdot{\bf h}}{{\bf |h_0|}}
\;.\nonumber\eea
The function $E$ lies in $L^1(\gO)$, since each individual term does.

We estimate the $L^1$-norm of $E$ as follows: let

\be U=\bigg\{x\in\gO\:\Big|\: \max\{|f(x)|,|g(x)|,|{\bf h}(x)|\} \geq m \bigg\}\;. \nonumber\ee
On $\gO\backslash U$, we can apply Taylor's Theorem to $F$ to find that

\be |E(x)|\leq C\left(f(x)^2+g(x)^2+|{\bf h}(x)|^2\right)\;,\ee
where $C$ depends on an upper bound for the second order partial derivatives in (\ref{partials}).
As a result,

\bea \int_{\gO\backslash U}|E(x)|\:dx \leq C\int_{\gO\backslash U} \Big(f(x)^2
+g(x)^2+|{\bf h}(x)|^2\Big)\:dx \nonumber\\
\leq C\bigg(\mathrm{Vol}(\gO)\norm{f}_{L^\infty(\gO)}^2+\norm{g}_{L^2(\gO)}^2
+\norm{{\bf h}}_{L^2(\gO)^n}^2\bigg) \;.\nonumber\eea
On $U$, we have

\bea \int_U |E|\leq C\int_U \bigg{(}|g_0+g|^2+|{\bf h_0}+{\bf h}|^2 \nonumber\\
+|F(f_0,g_0,|{\bf h_0}|)|+|\nabla F(f_0,g_0,|{\bf h_0}|)|(|f|+|g|+|{\bf h}|)\bigg{)} \;.\eea
The constant in front depends only on $F$ and $\norm{f_0}_{L^\infty(\gO)}$, as we need to consider only
$\norm{f}_{L^\infty(\gO)}\leq m$, say, which makes the constant $C(y)$ in (\ref{cofy}) bounded.
Using the fact that $f_0$, $g_0$, and ${\bf h_0}$ are bounded functions, we have

\bea \int_U |E|\leq C\int_U \left(1+|f|+|g|+|g|^2+|{\bf h}|+|{\bf h}|^2\right) \nonumber\\
\leq C\Big(\left(1+\norm{f}_{L^\infty(\gO)}\right)\mathrm{Vol}(U)\nonumber\\
+\left(\norm{g}_{L^2(\gO)}
+\norm{{\bf h}}_{L^2(\gO)^n}\right)\mathrm{Vol}(U)^{1/2} \nonumber\\
+\norm{g}_{L^2(\gO)}^2+\norm{{\bf h}}_{L^2(\gO)^n}^2\Big) \eea
Owing to the inequality

\be \mathrm{Vol}(U)\leq \frac{2}{m^2}\left(\norm{g}_{L^2(\gO)}^2 +\norm{{\bf h}}_{L^2(\gO)^n}^2\right) \nonumber\ee
for $\norm{f}_{L^\infty(\gO)}<m$,
we have

\be \int_U |E(x)|\:dx\leq C\left(\norm{f}_{L^\infty(\gO)}+\norm{g}_{L^2(\gO)}^2 +\norm{{\bf h}}_{L^2(\gO)^n}^2\right) \;.\ee
This proves the Fr\'{e}chet differentiability of the map $(f,g,{\bf h})\mapsto F(f,g,|{\bf h}|)$ at $(f_0,g_0,{\bf h_0})$.
Hence by the chain rule for Fr\'{e}chet derivatives, $F(\sigma, u, \abs{\nabla u})$ is Fr\'{e}chet differentiable
as a function of $\sigma$ at $\sigma_0$.
$\square$

\subsection{Proof of Lemma \ref{L:diffF3}}
{\bf Proof.}
The Fr\'echet differentiability of $F_{11}$ and $F_{22}$ follows from Lemma \ref{L:diff_general}.
Only Fr\'{e}chet differentiability of $F_{12}$ needs to be proven.

As a result of Lemma \ref{frdiff}, the map

\be \sigma\mapsto \left(\begin{array}{c}\nabla u^{(1)}_0 \\ \nabla u^{(2)}_0 \end{array} \right) \nonumber\ee
is Fr\'{e}chet differentiable from $L^\infty_\mathrm{ad}(\gO)$ to $L^2(\gO)^{2n}$ at $\sigma_0$.
We claim that the map

\bea L^2(\gO)^{2n}&\rightarrow& L^1(\gO) \nonumber\\
\label{dotproduct} \left(\begin{array}{c}v_1 \\ v_2 \end{array} \right) &\mapsto& \phi|v_1\cdot v_2|^{p/2} \eea
is Fr\'{e}chet differentiable at a pair of smooth vector fields $v_1$ and $v_2$ that satisfy

\be \label{bigproduct}
|v_1|,\; |v_2| \leq M,
 |v_1\cdot v_2| \geq \ga>0
\ee
on the support of $\phi$ for some $M>1$.
To show this, let $w_1$, $w_2\in L^2(\gO)^n$ and define a function $E\in L^1(\gO)$ by

\bea  E&=&\phi\bigg(\left|(v_1+w_1)\cdot(v_2+w_2)\right|^{p/2}-|v_1\cdot v_2|^{p/2}  \nonumber\\
\label{E:tayexp}
&-&\frac{p}{2}|v_1\cdot v_2|^{p/2-2}(v_1\cdot v_2)
(v_2\cdot w_1+v_1\cdot w_2)\bigg) \;.\eea

In order to estimate the $L^1(\gO)$-norm of $E$, define a set $U$ by

\be U=\left\{ x\in\gO\:\big|\: \max_{i=1,2} |w_i(x)|\geq \frac{\ga}{4M} \right\}\;.\nonumber\ee
On $\gO\backslash U$, both $w_1$ and $w_2$ are bounded above by $\ga/4M$, so we have

\bea |(v_1+w_1)\cdot (v_2+w_2)|&\geq& \ga-|v_1\cdot w_2|-|v_2\cdot w_1|-|w_1\cdot w_2| \nonumber\\
&\geq& \ga-|v_1|\frac{\ga}{4M}-|v_2|\frac{\ga}{4M}-\left(\frac{\ga}{4M}\right)^2 \nonumber\\
&\geq& \frac{\ga}{4} \;.\eea
It therefore follows from (\ref{E:tayexp}) and Taylor's formula applied to the function of $2n$ variables $|x\cdot y|^{p/2}$
(which is smooth and has bounded derivatives when the arguments satisfy $\ga/4\leq x,y\leq M$) that

\be \int_{\gO\backslash U} |E(x)|\:dx \leq C\left(\norm{w_1}_{L^2(\gO)^n}^2+\norm{w_2}_{L^2(\gO)^n}^2\right)\;, \ee
where $C$ depends on $m$, $M$, $p$, $n$, and $\ga$.

On $U$, we use the triangle inequality:

\bea |E|&\leq&|v_1\cdot v_2|^{p/2}\bigg|\Big|1+\frac{v_1\cdot w_2+v_2\cdot w_1+w_1\cdot w_2}{v_1\cdot v_2}\Big|^{p/2}-1\bigg|
\nonumber\\
\label{eofx}&+&\frac{p}{2}|v_1\cdot v_2|^{p/2-1}|v_2\cdot w_1+v_1\cdot w_2| \;.\eea
We make use of the inequality

\be \label{cp}
\left||1+z|^{p/2}-1\right|\leq C_p|z| \;.\ee
to estimate the term in the first line of (\ref{eofx}).
(For $|z|\leq1/2$ this inequality follows from the fact that the function $z\mapsto |1+z|^{p/2}$ has bounded derivative for
$|z|\leq1/2$, while for $|z|>1/2$ it follows easily from the fact that $p/2\leq1$.)
After using (\ref{cp}),
we obtain

\bea \int_U |E(x)|\:dx &\leq&\int_U C_p|v_1\cdot v_2|^{p/2-1}|v_1\cdot w_2+v_2\cdot w_1+w_1\cdot w_2| \nonumber\\
&+&\frac{p}{2}|v_1\cdot v_2|^{p/2-1}|v_2\cdot w_1+v_1\cdot w_2|\;. \eea
By the Cauchy-Schwarz inequality and the bounds on $v_1$, $v_2$, and $v_1\cdot v_2$,

\bea  \int_U |E(x)|\:dx &\leq& C\int_U |w_1|+|w_2|+|w_1\cdot w_2| \nonumber\\
&\leq& C\big(\norm{w_1}_{L^2(U)^n}\mathrm{Vol}(U)^{1/2}+\norm{w_2}_{L^2(U)^n}\mathrm{Vol}(U)^{1/2}\nonumber\\
&+&\norm{w_1}_{L^2(U)^n}\norm{w_2}_{L^2(U)^n}\big)\;. \eea

The volume of $U$ can be estimated by

\bea \mathrm{Vol}(U)&\leq& \mathrm{Vol}\Big(\left\{|w_1|\geq \frac{\ga}{4M}\right\}\Big)
+\mathrm{Vol}\Big(\left\{|w_2|\geq \frac{\ga}{4M}\right\}\Big) \nonumber\\
&\leq& 2\left(\frac{4M}{\ga}\right)^2\left(\norm{w_1}_{L^2(\gO)^n}^2+\norm{w_2}_{L^2(\gO)^n}^2\right) \;.\eea
Hence,

\be \int_U |E(x)|\:dx\leq C\left(\norm{w_1}_{L^2(\gO)^n}^2+\norm{w_2}_{L^2(\gO)^n}^2\right) \ee
This proves that the map (\ref{dotproduct}) is Fr\'{e}chet differentiable.  Since condition (\ref{bigproduct})
holds for the vector fields $\nabla u_0^{(1)}$ and $\nabla u_0^{(2)}$
on the support of $\phi$, the Fr\'{e}chet differentiability of $F_{12}$
at $\sigma_0$ follows from the chain rule.
$\square$

\subsection{Left parametrix for over-determined DN elliptic operator}
We provide here a sketch of the classical construction of the left parametrix $B$ used in the proof of Theorem \ref{T:qpat}.

According to the assumptions of Theorem \ref{T:qpat}, there exists a finite open covering $\{U_j\}$ of the compact subset $\overline{\Omega'}\times S^{n-1}$ in $\Omega''\times \left(\R^{n}\setminus \{0\}\right)$, such that for each $j$, there exists a number $k_j$ such that $A_{k_j} (x,\xi)$ is invertible for $(x,\xi)\in U_j$. Here $S^{n-1}$ is the unit sphere in the $\xi$-space $\R^n$.  Consider a smooth partition of unity $\psi_j(x,\xi)$ on $\overline{\Omega'}\times S^{n-1}$ subordinated to the covering. We can always assume that it is positively homogeneous of order zero with respect to $\xi$ outside a neighborhood of the origin and smooth on $\Omega''\times \R^{n}$. Let us also denote by $P_k:\C^{2n}\mapsto \C^2$ the operator such that $P_k(a_1,a_2,\dots,a_{2n-1},a_{2n})=(a_{2k-1},a_{2k})$. One sees that $A_k=P_kA$. Then one can check that the symbol $B(x,\xi)=\sum_j \psi_j (x,\xi) A_{k_j}(x,\xi)^{-1} P_{k_j}$ is the one we require.

\section{Final remarks}\label{S:remarks}
\begin{enumerate}
\item As we have mentioned, some of the models discussed in this article have been studied previously.  Reconstructions from the functionals (\ref{E:2data}) with $p=2$
have been performed for instance in \cite{bal3,cap,KuKuAET}, with similar stability estimates being obtained in \cite{bal3,KuKuAET}, the approach
in \cite{bal3} being global.
Global reconstruction from a single functional was considered in \cite{bal1} by solving a Cauchy problem inward from parts of $\dgO$ for a nonlinear hyperbolic equation.
In this case, stability results were obtained in two dimensions, and in parts of $\gO$ in higher dimensions.
\item The $p=1$ case has been studied in \cite{Nach1,Nach2,Nach3}.  In \cite{Nach1},
an iterative reconstruction procedure was provided, whose effectiveness was demonstrated in numerical experiments.  Our analysis in terms of pseudo-differential operators
implies that although the inversion of the corresponding linearized operator from a single functional is not an elliptic problem, ellipticity only fails at points $(x,\xi)$ where
$\nabla u_0\parallel\xi$.  Since most singularities of $\rho$ would likely not satisfy that condition, we should expect accurate reconstruction of sharp features
in almost all cases, as the numerical experiments demonstrate.
\item While we have shown infinitesimal Fredholm property of the problems with internal data, infinitesimal uniqueness has not been shown and we suspect that it does not hold under our very general conditions. However, uniqueness should hold {\bf generically}, which we plan to address elsewhere.
\item Besides absence of an infinitesimal uniqueness result, there is another obstacle for obtaining the local uniqueness and stability result for the non-linear problem. Namely, Fr\'{e}chet differentiability is proven in worse function spaces than the Fredholm property. We also plan to address this discrepancy in a future work.
\item We suggest that parametrices constructed in this paper could be used for approximate reconstructions and for pre-conditioning iterative methods.
\item In the cases when our ellipticity analysis asks for multiple measurements, this does not mean that reconstruction with a smaller number of measurements is impossible.
On the contrary, such reconstructions have been achieved by solving hyperbolic and degenerate elliptic problems in \cite{bal2,KuKuSynt,KuKuAET}.
However, such approaches naturally lead to correspondingly error propagation and some blurring effects for the parts of the wavefront sets where ellipticity is lost \cite{KuKuAET,Kuch_hybrid}.
\item The types of the internal data functionals $F$ considered in this paper do not cover all the needs of hybrid methods.
For instance, some non-local functionals of $\sigma$ arise in UOT \cite{AllmBang}. We plan to address those in a subsequent work.
\end{enumerate}

\section*{Acknowledgments} The work of P.~K. was partially supported by the NSF DMS Grant \# 0604778. The work of both authors was supported in part by the Award No. KUS-C1-016-04, made by King Abdullah University of Science and Technology (KAUST) and by the IAMCS.
The authors also wish to thank the referees for their very helpful suggestions and remarks.

\end{document}